\documentclass[psamsfonts, reqno]{amsart}

\usepackage{amscd,amsmath,amssymb,amsfonts,verbatim}
\usepackage{times}
\usepackage[cmtip, all]{xy}
\usepackage{graphicx}
\usepackage[active]{srcltx}
\usepackage{color}
\usepackage{textcomp}

\definecolor{refkey}{gray}{.5}   
\definecolor{labelkey}{gray}{.5} 
\definecolor{Red}{rgb}{1,0,0}

\newcommand{\pf}{{\bf Proof : }}

\newcommand{\qedwhite}{\hfill \ensuremath{\Box}} 

\newtheorem{theo}{Theorem}[section]

\newtheorem{lem}[theo]{Lemma}
\newtheorem{cor}[theo]{Corollary}

\theoremstyle{definition}

\newtheorem{notn}[theo]{Notation}
\newtheorem{rem}[theo]{Remark}
\newtheorem{exam}[theo]{Example}
\newtheorem{defi}[theo]{Definition}

\title{On the first row map ${GL}_{d+1}(R)\longrightarrow \frac{{Um}_{d+1}(R)}{{E}_{d+1}(R)}$}
\author{Sampat Sharma}
\newcommand{\Addresses}{{
  \bigskip
  \footnotesize

   \textsc{Sampat Sharma, Department of Mathematics, IIT Bombay,\\ \noindent
           Main Gate Rd, IIT Area, Powai, Mumbai, Maharashtra 400076
   } \par\nopagebreak
  \textit{E-mail}: Sampat ~Sharma \texttt{<sampat@math.iitb.ac.in>; <sampat.iiserm@gmail.com>}

  \medskip

  }}

\begin{document}
\maketitle
\subjclass 2020 Mathematics Subject Classification:{13H99, 13E05, 19B99}

 ~ {Unimodular rows, Group homomorphism}
 \begin{abstract}
 In this paper, we prove that there is a group homomorphism from general linear group over a polynomial extension 
 of a local ring to the W. van der Kallen's group $\frac{Um_{d+1}(R[X])}{E_{d+1}(R[X])}.$
 \end{abstract}
 
 \vskip0.50in

\section{Introduction} 
 Throughout this article, we assume $R$ to be a commutative noetherian ring with $1 \neq 0 $ unless stated otherwise. 
\par Group homomorphisms from special linear group (or from general linear group) to orbit sets of unimodular rows plays an important role in solving problems in classical $K$-theory. In $($\cite[Theorem 3.16(iv)]{vdk1}$),$ W. van der Kallen proves that for a commutative ring $R$ of dimension $d$, the 
first row map 
$$SL_{d+1}(R)\longrightarrow \frac{Um_{d+1}(R)}{E_{d+1}(R)}$$ 
is a group homomorphism. This result was crucial for Ravi A. Rao to settle the Bass--Quillen conjecture in dimension 3 in \cite{invent}.

Further, examples has been given by W. van der Kallen showing that the first row map is not a group homomorphism in general. To be precise, for a commutative ring $R$ of dimension $d$ the first row map 
$$GL_{d+1}(R)\longrightarrow \frac{Um_{d+1}(R)}{E_{d+1}(R)}$$ 
is not a group homomorphism $($\cite[Examples 4.16, 4.13]{vdk1}$).$ W. van der Kallen also gave an example in $(${\cite [Proposition 7.10]{vdk2}}$)$ where he shows that the first 
row map 
$$GL_{d}(R)\longrightarrow \frac{Um_{d}(R)}{E_{d}(R)}$$ 
is not a group homomorphism for a commutative ring $R$ of dimension $d.$

\par In this article, we show that the first row map at $d+1$ level is a group homomorphism when we work over 
polynomial extension 
of a local ring. In particular, 
we prove the following result: 

 \begin{theo}
 \label{maintheorem}
 Let $R$ be a reduced local ring of dimension $d, d\geq 2, \frac{1}{d!}\in R.$ Let $v\in Um_{d+1}(R[X])$ and $g\in GL_{d+1}(R[X])$, 
 then $[vg] = [v]\ast [g].$ 
 Consequently, the first row map 
 $$GL_{d+1}(R[X]) \longrightarrow \frac{Um_{d+1}(R[X])}{E_{d+1}(R[X])}$$  
 is a group homomorphism.
\end{theo}

\begin{cor}\label{coro}
   Let $R$ be a reduced local ring of dimension $d, d\geq 2, \frac{1}{d!}\in R.$ Then the set 
$GL_{d}(R[X])E_{d+1}(R[X])$ contains the commutator subgroup of $GL_{d+1}(R[X]).$
\end{cor}

\par It is interesting to note here that, if instead we work over a ring $R$ of dimension $d$, the commutator subgroup $[GL_{d+1}(R), GL_{d+1}(R)]$ need not to be contained in $GL_{d}(R)E_{d+1}(R)$ as shown by W. van der Kallen in 
$(${\cite[Example 4.16]{vdk1}}$).$

\section{Preliminaries}
\par 
Let $v = (a_{0},a_{1},\ldots, a_{r}), w = (b_{0},b_{1},\ldots, b_{r})$ be two rows of length $r+1$ over a commutative ring 
$R$. A row $v\in R^{r+1}$ is said to be unimodular if there is a $w \in R^{r+1}$ with
$\langle v ,w\rangle = \Sigma_{i = 0}^{r} a_{i}b_{i} = 1$ and $Um_{r+1}(R)$ will denote
the set of unimodular rows (over $R$) of length 
$r+1$.
\par 
The group of elementary matrices,  denoted by $E_{r+1}(R)$,  is a subgroup of $GL_{r+1}(R)$ and is generated by the matrices 
of the form $e_{ij}(\lambda) = I_{r+1} + \lambda E_{ij}$, where $\lambda \in R, ~i\neq j, ~1\leq i,j\leq r+1,~
E_{ij} \in M_{r+1}(R)$ whose $ij^{th}$ entry is $1$ and all other entries are zero. The elementary linear group 
$E_{r+1}(R)$ acts on the rows of length $r+1$ by right multiplication. Moreover, this action takes unimodular rows to
unimodular
rows and we denote set of orbits of this action by $\frac{Um_{r+1}(R)}{E_{r+1}(R)}.$ 
The equivalence class of a row $v$ under this equivalence relation is denoted by $[v].$

\begin{defi}{\bf{Bass--Serre dimension:}}
 Let $R$ be a ring whose maximal spectrum $\mbox{Max}(R)$ is a finite union of subsets $V_{i},$ where each $V_{i},$ 
 when endowed with the (topology induced from the) Zariski topology, is a space of Krull dimension $d.$ We shall say Bass--Serre dimension of 
 $R$ is $d$ in such a case.
\end{defi}

\begin{exam}
 Let $R$ be a reduced local ring of dimension $d\geq 1,$ then Bass--Serre dimension of $R[X]$ is $d$ as for any non-unit,
 non-zero-divisor 
 $\pi \in R$ 
 $$\mbox{Spec}(R[X]) = \mbox{Spec}(\frac{R}{(\pi)}[X]) \cup \mbox{Spec}(R_{\pi}[X])$$ 
 is a finite union of noetherian spaces of dimension $d.$
 \end{exam}

\par 
In $(${\cite[Theorem 3.6]{vdk1}}$),$ W. van der Kallen derives an abelian
 group structure on $\frac{Um_{d+1}(R)}{E_{d+1}(R)}$ when Bass--Serre dimension of
$R$ is $d,$ for all $d\geq 2.$ We will denote the group operation in this group by $\ast.$

We recall a few definitions to be used in the next section.
\begin{defi}
 Let $I$ be an ideal of a polynomial ring $R[X]$ (in one indeterminate). By $l(I)$ we denote the set consisting of $0$ and 
 all leading coefficients of $f\in I\diagdown \{0\}.$ Obviously $l(I)$ is an ideal of $R.$
\end{defi}

\begin{defi}
For a prime ideal $\mathfrak{P}$ of a commutative ring $R$, the height of $\mathfrak{P}$, denoted by $\mbox{ht}_{R}{\mathfrak{P}}$, is defined by 
$$\mbox{ht}_{R}{\mathfrak{P}} = \mbox{sup}\{r ~| ~\mbox{there}~\mbox{exists}~\mbox{a}~\mbox{chain}~\mathfrak{P_{0}} \subsetneq \mathfrak{P_{1}}\subsetneq \cdots \subsetneq \mathfrak{P}_{r} = \mathfrak{P}~\mbox{in}~\mbox{Spec}~R\}.$$
\end{defi}

\begin{defi} The height of a proper ideal $I$ of a commutative ring $R$, denoted by $\mbox{ht}_{R}I$, is defined by 
$$\mbox{ht}_{R}{I} = \mbox{inf}\{\mbox{ht}_{R}{\mathfrak{P}}~ |~ \mathfrak{P}\in \mbox{Spec}~R, I\subset \mathfrak{P}\}.$$
\end{defi}

\begin{rem}In a commutative noetherian ring $R$ there are only finitely many prime ideals minimal over a given ideal $I$. It follows that 
$$\mbox{ht}_{R}{I} = \mbox{min}\{\mbox{ht}_{R}{\mathfrak{P}}~ | ~\mathfrak{P} ~\mbox{minimal}~\mbox{over}~I\}.$$
\end{rem}

We recall the following lemma by H. Bass and A. Suslin $(${\cite[Chapter 3, Lemma 3.2]{lamoldbook}}$).$
\begin{lem}
\label{heightleading}
 Let $R$ be a commutative noetherian ring and $I$ be an ideal of $R[X].$ Then $\mbox{ht}_{R}l(I)\geq \mbox{ht}_{R[X]}I.$
\end{lem}

\begin{lem}
 \label{makingnzd} Let $R$ be a reduced noetherian ring and $ f = (f_{0}, f_{1}, \ldots, f_{r})\in Um_{r+1}(R[X]), r\geq 1.$ Then there 
 exists $g = (g_{0}, g_{1}, \ldots, g_{r})\in Um_{r+1}(R[X])$ in the elementary orbit of $f$ such that $l(g_{0}) = \pi,$ 
  a non-zero-divisor in $R.$
\end{lem}
${\pf}$ By $(${\cite[Corollary 9.4]{7}}$),$ there exists $\varepsilon \in E_{r+1}(R[X])$ such that 
$(f_{0}, \ldots, f_{r})\varepsilon = (g_{0}, \ldots, g_{r})$ and $\mbox{ht}(g_{1}, \ldots, g_{r}) \geq r \geq 1.$ By 
Lemma \ref{heightleading}, $\mbox{ht}_{R}l(I)\geq \mbox{ht}_{R[X]}I \geq r \geq 1,$ where $I = \langle g_{1}, \ldots, g_{r}\rangle.$ In a reduced noetherian ring, the set of zero divisors is the union of minimal prime ideals. Thus 
by the prime avoidance lemma 
we observe that the ideal $l(I)$ contains a non-zero-divisor.
Therefore, there exist $\lambda_{i} \in R$ such that $\pi = \sum_{i=1}^{r}\lambda_{i}l(g_{i})$ is a non-zero-divisor in $R.$ 
 Let $\mbox{deg}(g_{i}) = d_{i}$ for $0\leq i\leq r.$ Replace $g_{0}$ with $g_{0} + \sum_{i=1}^{r} \lambda_{i}X^{m-d_{i}+1}g_{i}$, where $m \geq \{d_{i}\}_{i\geq 0},$ to get $l(g_{0}) = \pi.$
$~~~~~~~~~~~~~~~~~~~~~~~~~~~~~~~~~~~~~~~~~~~~~~~~~~~~~~~~~~~~~~~~~~~~~~~~~~~~~~~~~~~~~~~~~~~~~~~~~~~~~~~~~~~~~~\qedwhite$

The following three results are due to Roitman $($\cite[Lemma 1, Lemma2, Lemma 3]{roitstably}$).$ 
\begin{lem}
 \label{roit1} Let $(x_{0}, x_{1}, \ldots, x_{r})\in Um_{r+1}(R),$ $r\geq 2$ and $t$ be an element of $R$ which is invertible 
 modulo $(x_{0}, \ldots, x_{r-2})$. Then, 
 $$(x_{0}, \ldots, x_{r}) \underset{E}{\sim} (x_{0}, \ldots, x_{r-1}, t^{2}x_{r}).$$
\end{lem}

\begin{lem}
 \label{roit2} Let $S$ be a multiplicative subset of $R$ such that $R_{S}$ is a noetherian ring of finite Krull dimension $d.$ Let 
 $(\overline{a_{0}}, \ldots, \overline{a_{r}}) \in Um_{r+1}(R_{S})$, $r>d.$ Then there exist $b_{i}\in R$ $(1\leq i\leq r)
 $ and $s\in S$ such that $s\in R(a_{1}+b_{1}a_{0}) + \cdots + R(a_{r} + b_{r}a_{0}).$
 \end{lem}
 
 \begin{lem}
 \label{roit3}
  Let $f(X) \in R[X]$ have degree $n>0,$ and let $f(0)$ be a unit. Then for any $g(X) \in R[X]$ and any natural 
 number $k\geq (\mbox{deg}(g(X))-\mbox{deg}(f(X)) + 1)$, there exists $h(X) \in R[X]$ of degree $< n$ such that 
 $$g(X) = X^{k}h(X)~\mbox{modulo}~(f(X)).$$
\end{lem}

The next result is due to Suslin and Vaserstein $($\cite[Lemma 11.1]{7}$).$ 
\begin{lem}
 \label{susmonic} Let $I$ be an ideal of $R[X]$ containing a monic polynomial $f$ of degree $m$. Let $g_{1}, \ldots, g_{k} \in I
 $ with $\mbox{deg}(g_{i}) < m,$ for $1\leq i\leq k.$ Assume that the coefficients of the $g_{i},$ $1\leq i\leq k,$ generate $R$. 
 Then $I$ contains a monic polynomial of degree $m-1.$
\end{lem}

\begin{lem}\label{susmoncor} Let $I$ be an ideal of $R[X]$ containing a monic polynomial. Let $m\geq 1.$ Let $g_{1}, \ldots, g_{k} \in I
 $ with $\mbox{deg}(g_{i}) < m,$ for $1\leq i\leq k.$ Assume that the coefficients of the $g_{i},$ $1\leq i\leq k,$ generate $R$. 
 Then $I$ contains a monic polynomial of every degree $\geq m-1.$
\end{lem}
${\pf}$ By Lemma \ref{susmonic}, if $I$ contains a monic polynomial of degree $n$ with $n\geq m,$ then it contains one of degree $n-1.$ Thus $I$ contains a monic polynomial of every degree $\geq m-1.$ 
$~~~~~~~~~~~~~~~~~~~~~~~~~~~~~~~~~~~~~~~~~~~~~~~~~~~~~~~~~~~~~~~~~~~~~~~~~~~~~~~~~~~~~~~~~~~~~~~~~~~~~~~~~~~~~~\qedwhite$

\section{The main result}

\begin{defi}
A polynomial $f(X)\in R[X]$ is said to be a $\pi$-power monic polynomial if its highest coefficient is $\pi^{k}$ for some $k.$
\end{defi}

Let $v\in Um_{r+1}(R[X]),$ where $r\geq \frac{d}{2} + 1,$ with $R$ a local ring of dimension $d.$ M. Roitman's argument in 
$(${\cite[Theorem 5]{roitstably}}$)$ and Ravi Rao's 
argument in $(${\cite[Proposition 4.8]{twoexam}}$)$ show how one could decrease the 
degree of all but one (special) coordinate of $v.$ In $(${\cite[Proposition 1.4.4]{invent}}$)$ Ravi Rao suitably modified M. Roitman's argument and proved that for $r\geq \mbox{max}(2, d/2 + 1),$ unimodular rows of length $r+1$ 
over a polynomial extension of a local ring $R$ can be elementarily mapped to a unimodular row such that all except first two of its coordinates lie in $R$ and the last coordinate is a non-zero-divisor.
 In $(${\cite[Theorem 4.10]{sampat}}$)$, I have also modified M. Roitman's argument to prove the following result. 
\begin{theo}
 Let $R$ be a reduced finitely generated ring over $\mathbb{Z}$ of dimension $d, d\geq 3.$ Let  
 $u(X) = (u_{0}(X), u_{1}(X), \ldots, u_{d}(X))
  \in 
 Um_{d+1}(R[X]).$ Then 
 $$(u_{0}(X), u_{1}(X), \ldots, u_{d}(X))\overset{E_{d+1}(R[X])}{\sim} (v_{0}(X), v_{1}(X), c_{2}, \ldots, c_{d})$$
 for some $(v_{0}(X), v_{1}(X), c_{2}, \ldots, c_{d})\in Um_{d+1}(R[X])$, $c_{i}\in R, 2\leq i\leq d$, and 
 $c_{d}$ is a non-zero-divisor.
\end{theo}

The proof of the following Theorem is modeled after the proof of Roitman's claim on $(${\cite[p. 587]{roitstably}}$).$

\begin{theo}
 \label{roittype}
 Let $R$ be a reduced local ring of dimension $d, d\geq 2.$ Let $a$ be an invertible element in $R$ and 
 $u(X) = (au_{0}(X), u_{1}(X), \ldots, u_{d}(X))
  \in 
 Um_{d+1}(R[X]).$ Then 
 $$(au_{0}(X), u_{1}(X), \ldots, u_{d}(X))\overset{E_{d+1}(R[X])}{\sim} (av_{0}(X), v_{1}(X), c_{2}, \ldots, c_{d})$$
 for some $(v_{0}(X), v_{1}(X), c_{2}, \ldots, c_{d})\in Um_{d+1}(R[X])$, $c_{i}\in R, 2\leq i\leq d$, $\mbox{deg}(v_{1}(X) ) = 1$ and 
 $c_{d}$ is a non-zero-divisor.
\end{theo}
${\pf}$  In view of Lemma \ref{makingnzd}, we may assume that the leading coefficient $l(u_{0})$ of $u_{0}$ equals some non-zero-divisor $\pi.$ If $\pi$ is a unit, then upon using Whitehead lemma $($\cite[Corollary 2.3]{7}$)$ we can 
make $u_{0}(X)$ a monic polynomial. Therefore in view of Suslin's result $($\cite[Corollary 4.4.8]{ideals}$)$ we can elementarily map $u(X)$ to $(1,0,\ldots, 0).$ Thus we are done. 
\par Let us now assume that $\pi$ is a non-unit. Let
$\overline{R} = \frac{R}{\pi R}$ and we elementarily map $\overline{u(X)}$ to $(1, 0, \ldots, 0)$ with the help of $(${\cite[Theorem 7.2]{4}}$).$ Upon lifting the elementary map, we obtain a row 
$v(X) = (v_{0}(X), v_{1}(X), \ldots, v_{d}(X))$ by applying elementary transformations to $u(X)$ such that 
$$v(X) \equiv (1, 0, \ldots, 0)~\mbox{mod}(R[X]\pi).$$ 
\par We can perform such trasformations so that at every stage the row contains a $\pi$-power monic polynomial. Indeed, if for example we have to perform the elementary transformation 
$$(h_{0},h_{1}, \ldots, h_{d}) \overset{T}{\longrightarrow} (h_{0}, h_{1}+bh_{0}, h_{2}, \ldots, h_{d})$$ 
and $h_{1}$ is a $\pi$-power monic polynomial, then we replace $T$ by the following transformations 
\begin{align*}
 (h_{0}, h_{1}, \ldots, h_{d})
 & \longrightarrow (h_{0}+\pi X^{n}h_{1}, h_{1}, h_{2}, \ldots, h_{d})\\
 & \longrightarrow (h_{0}+ \pi X^{n}h_{1}, h_{1}+ b(h_{0}+\pi X^{m}h_{1}) , h_{2}, \ldots, h_{d})
\end{align*}
where $n > \mbox{deg}(h_{0}).$

Going back to the proof of the theorem we assume that 
$$(au_{0}, u_{1}, \ldots, u_{d}) \equiv (1, 0,  \ldots, 0)~\mbox{mod}(R[X]\pi)$$ 
and $u_{i}$ is a $\pi$-power monic polynomial. If $i> 0$ we replace $au_{0}$ with $au_{0} + \pi X^{m}u_{i}.$ So 
we assume that $u_{0}$ is a $\pi$-power monic polynomial and $\mbox{deg}(u_{0}) > 0.$ Since any unimodular row over a local ring can be elementarily mapped to $(1,0, \ldots, 0),$ we may also assume that $u_{0}(0) = 1.$ Now, choose 
$k$ such that $$2k \geq \mbox{max}\{\mbox{deg}(u_{i}(X)) - \mbox{deg}(u_{0}(X)) + 1, ~~~~1\leq i\leq d\}.$$
By Lemma \ref{roit3}, we can assume 
$u_{i} = X^{2k}h_{i},$ where $\mbox{deg}(h_{i}) < \mbox{deg}(u_{0})$ for $1\leq i\leq d.$ Since $X$ is invertible modulo $u_{0}(X)$, upon applying Lemma \ref{roit1} with $t = X^{-1}$, we may assume 
$\mbox{deg}(u_{i}) < \mbox{deg}(u_{0})$ for $1\leq i\leq d.$
\par Let $\mbox{deg}( u_{0}) = m.$ If $m = 1, $ then $u_{i}\in R$ for $1\leq i\leq d.$ Now suppose that $m\geq 2.$ 
Let $(c_{1}, c_{2}, \ldots, c_{m(d-1)})$ be the coefficients of $1, X, \ldots, X^{m-1}$ in the polynomials 
$u_{2}(X), \ldots, u_{d}(X).$ By $($\cite[Chapter III, Lemma 1.1]{lamoldbook}$),$ the ideal generated in $R_{\pi}$ by 
$R_{\pi} \bigcap (R[X]_{\pi}\overset{\sim}{u_{0}} + R[X]_{\pi}\overset{\sim}{u_{1}})$ and the coefficients of $\overset{\sim}{u_{i}}$ $(
2\leq i\leq r)$ is $R_{\pi}.$ Here $\overset{\sim}{u_{i}}$ denotes the image of $u_{i}$ in $R[X]_{\pi}.$ Since $d\geq 2,$ we have $m(d-1)> \mbox{dim}R_{\pi}.$ Thus, by Lemma \ref{roit2}, there exist 
$$(c_{1}', c_{2}', \ldots, c_{m(d-1)}') \equiv (c_{1}, c_{2}, \ldots, c_{m(d-1)}) ~ \mbox{mod} 
~ ((R[X]u_{0} + R[X]u_{1}) \bigcap R) $$ such that 
$R_{\pi}\overset{\sim}{c_{1}'} + \cdots + R_{\pi}\overset{\sim}{c_{m(d-1)}'} = R_{\pi}.$ Adding suitable multiples of $u_{0}$ and $u_{1}$ to $u_{i}$, $2\leq i\leq d$, we can reduce the situation to the case 
$R_{\pi}\overset{\sim}{c_{1}} + \cdots + R_{\pi}\overset{\sim}{c_{m(d-1)}} = R_{\pi}.$ By Lemma \ref{susmoncor}, the ideal 
$R_{\pi}[X]\overset{\sim}{u_{0}} + R_{\pi}[X]\overset{\sim}{u_{2}} + \cdots + R_{\pi}[X]\overset{\sim}{u_{d}}$ contains a monic polynomial of degree $m-1.$ This implies that $R[X]u_{0} + R[X]u_{2} + \cdots + R[X]u_{d}$ contains a 
 $\pi$-power monic polynomial $h(X)$ of degree $m-1.$ Suppose that the leading coefficient of $h(X)$ is ${\pi}^{k}$. Now $u_{1}$ may be replaced with a $\pi$-power monic polynomial of degree $m-1.$ 
This does not depend on the original degree of $u_{1}.$ First one does the elementary transformations to bring the degree of $u_{1}$ down to $m-2$ and then simply adds $h$ to $u_{1}.$ Indeed, in the case $\mbox{deg}(u_{1}) = m-1$ 
we perform the following elementary transformations :
\begin{align*}
 (u_{0}, u_{1}, \ldots, u_{d}) & \longrightarrow (u_{0}, {\pi}^{2k}u_{1}, \ldots, u_{d})\\
 & \longrightarrow (u_{0}, {\pi}^{2k}u_{1}+ (1-{\pi}^{k}l(u_{1}))h,u_{2}, \ldots, u_{d}).
 \end{align*}

\par Note that we still have $u_{0} \equiv ~1~\mbox{mod}~\pi R[X]$.  Moreover we have $\mbox{deg}(u_{0}) = m$ and  $\mbox{deg}(u_{1}) = m-{1}.$  We can reduce the degree of the polynomials $u_{i}$, $2\leq i\leq d$ by 
subtracting suitable multiples of $u_{1}$ from $u_{i}$'s. In order to accomplish this, we may need at first to multiply each $u_{i}$ by a suitable power of $\pi$ by virtue of Lemma \ref{roit1}. Thus, we can assume that $\mbox{deg}(u_{i}) < m-1$ 
for $2\leq i\leq d.$
\par If $m-1 = 1, $ then we are done. Assume that $m-1\geq 2.$ Let $(t_{1}, t_{2}, \ldots, t_{(m-1)(d-1)})$ be the coefficients of $1, X, \ldots, X^{m-2}$ in the polynomials 
$u_{2}(X), \ldots, u_{d}(X).$ Repeating the argument above and using lemma \ref{susmoncor}, we can find a $\pi$-power monic polynomial $h_{1}(X)$ of degree $m-2$ in the ideal $R[X]u_{0} + R[X]u_{2} + \cdots + R[X]u_{d}.$ Now as above $u_{1}$ may be replaced with 
a $\pi$-power monic polynomial of degree $m-2.$ By Lemma \ref{roit1}, we also assume that 
$\mbox{deg}( u_{i}) < m-{2}$ for $2\leq i\leq d.$ Repetition of this argument yields $\mbox{deg}(u_{1}) = 1,$ $ u_{i} \in R$ for $2\leq i\leq d.$ Thus, it proves that $$(au_{0}(X), u_{1}(X), \ldots, u_{d}(X))\overset{E_{d+1}(R[X])}{\sim} (v_{0}(X), v_{1}'(X), c_{2}, \ldots, c_{d})$$
 for some $(v_{0}(X), v_{1}'(X), c_{2}, \ldots, c_{d})\in Um_{d+1}(R[X])$, $c_{i}\in R$ and $\mbox{deg}(v_{1}'(X)) = 1$ for  $2\leq i\leq d.$ Now using Whitehead lemma $($\cite[Corollary 2.3]{7}$),$ we have 
$$(au_{0}(X), u_{1}(X), \ldots, u_{d}(X))\overset{E_{d+1}(R[X])}{\sim} (av_{0}(X), v_{1}(X), c_{2}, \ldots, c_{d}),$$ 
where $v_{1}(X) = a^{-1}v_{1}'(X).$

\par We are only left to show that $c_{d}$ is a non-zero-divisor. Let us assume that $0\neq c_{d}$ is a zero divisor, then there exists $0\neq r\in R$ such that $c_{d}\cdot r = 0.$ Since $(av_{0}(X), v_{1}(X), c_{2}, \ldots, c_{d})\in Um_{d+1}(R[X])$, $ r \in \langle av_{0}(X), v_{1}(X), \ldots, c_{d-1}\rangle.$ Therefore, 
$$(av_{0}(X), v_{1}(X), c_{2}, \ldots, c_{d}) \overset{E_{d+1}(R[X])}{\sim} (av_{0}(X), v_{1}(X), c_{2}, \ldots,c_{d-1}, c_{d}+ r).$$ Since $R$ is a reduced noetherian ring, the set of zero divisors of $R$ is the union of minimal prime ideals of $R$. Let $\mathfrak{P_{1}}, \ldots, \mathfrak{P}_{n}$ be the minimal prime ideals of $R$. Suppose $c_{d}$ belongs to the first $l$ prime ideals (reenumerating if necessary), $l<n.$ Then $r\in \cap_{i=l+1}^{n}\mathfrak{P_{i}}$ and $r$ does not belong to at least one of the first $l$ prime ideals. Otherwise $r$ would be zero as $R$ is a reduced ring and in a reduced noetherian ring intersection of all minimal prime ideals is $\mbox{nil}(R).$ Therefore $c_{d} + r$ belongs to less than $l$ prime ideals. Now inducting on the number of minimal prime ideals containing $c_{d},$ we can make the last coordinate a non-zero-divisor.
$~~~~~~~~~~~~~~~~~~~~~~~~~~~~~~~~~~~~~~~~~~~~~~~~~~~~~~~~~~~~~~~~~~~~~~~~~~~~~~~~~~~~~~~~~~~~~~~~~~~~~~~~~
           ~~~~~\qedwhite$

\begin{lem}
 \label{unitequivalence}
 Let $R$ be a commutative ring and $v = (v_{1}, \ldots, v_{n}), ~ w = (w_{1}, \ldots, w_{n}) \in Um_{n}(R)$, $n\geq 3$ and 
 $u\in R^{\times}$
 be such that $ (uv_{1},v_{2}, \ldots, v_{n})\overset{E_{n}(R)}{\sim} (uw_{1},w_{2}, \ldots, w_{n}),$ then 
 $$(v_{1},v_{2},  \ldots, v_{n}) \overset{E_{n}(R)}{\sim} (w_{1}, w_{2},\ldots, w_{n}).$$
\end{lem}
${\pf}$ Let $\varepsilon \in E_{n}(R)$ be such that $ (uv_{1},v_{2}, \ldots, v_{n})\varepsilon = (uw_{1},w_{2}, \ldots, w_{n}).$ 
Therefore, $$(v_{1}, \ldots, v_{n})\begin{bmatrix}
                         u & 0\\
                         0 & I_{n-1}\\
                        \end{bmatrix}\varepsilon = (w_{1}, \ldots, w_{n})\begin{bmatrix}
                         u & 0\\
                         0 & I_{n-1}\\
                        \end{bmatrix}.$$
In view of $(${\cite[Corollary 1.4]{4}}$),$ 
$E_{n}(R) \trianglelefteq GL_{n}(R)$ for $n\geq 3,$ we have $v\varepsilon_{1} = w$ for some $\varepsilon_{1}\in E_{n}(R).$

     $~~~~~~~~~~~~~~~~~~~~~~~~~~~~~~~~~~~~~~~~~~~~~~~~~~~~~~~~~~~~~~~~~~~~~~~~~~~~~~~~~~~~~~~~~~~~~~~~~~~~~~~~~~~~~~\qedwhite$
     
     \begin{lem}
      \label{removingunit}
      Let $R$ be a reduced local ring of dimension $d, d\geq 2, \frac{1}{d!}\in R.$ Let $u\in R^{\times}$ and 
 $v(X) = (v_{0}(X), v_{1}(X), \ldots, v_{d}(X))
  \in 
 Um_{d+1}(R[X]).$ Then, 
 $$(uv_{0}(X), v_{1}(X), \ldots, v_{d}(X))\overset{E_{d+1}(R[X])}{\sim} (v_{0}(X), v_{1}(X), \ldots, v_{d}(X)).$$
      
     \end{lem}
${\pf}$ We will prove the lemma by induction on $d.$ Let $d=2$, $ v(X) = (v_{0}(X), v_{1}(X), v_{2}(X))\in Um_{3}(R[X]).$ 
By $(${\cite[Theorem 2.4]{invent}}$),$ $v(X)$ is completable. Thus, we are done 
 in view of $(${\cite[Lemma 1.5.1]{invent}}$).$  
 \par Let $d>2$ and $v(X) = (v_{0}(X), v_{1}(X), \ldots, v_{d}(X))
  \in 
 Um_{d+1}(R[X]).$ By Theorem \ref{roittype}, we have 
 \begin{equation}
 \label{eq1}
  (uv_{0}(X), v_{1}(X), \ldots, v_{d}(X))\overset{E_{d+1}(R[X])}{\sim} (uu_{0}(X), u_{1}(X), c_{2}, \ldots, c_{d})
 \end{equation}
 for some $(u_{0}(X), u_{1}(X), c_{2}, \ldots, c_{d})\in Um_{d+1}(R[X])$, $c_{i}\in R, 2\leq i\leq d$, and 
 $c_{d}$ is a non-zero-divisor. By Lemma \ref{unitequivalence}, we have 
 \begin{equation}
 \label{eq2}
  (v_{0}(X), v_{1}(X), \ldots, v_{d}(X))\overset{E_{d+1}(R[X])}{\sim} (u_{0}(X), u_{1}(X), c_{2}, \ldots, c_{d}).
 \end{equation}
 Let $\overline{R} = \frac{R}{c_{d}R}.$ Since $\mbox{dim}{\overline{R}} = d-1,$ by induction hypothesis one has 
  $$(\overline{uu_{0}(X)}, \overline{u_{1}(X)}, \overline{c_{2}}, \ldots, \overline{c_{d-1}})
\overset{E_{d}(\overset{-}{R}[X])}{\sim} 
(\overline{u_{0}(X)}, \overline{u_{1}(X)}, \overline{c_{2}}, \ldots, \overline{c_{d-1}}).$$ 
Lifting the elementary map and making appropriate elementary transformations, we obtain that 
\begin{equation}
 \label{eq3}
 (uu_{0}(X), u_{1}(X), c_{2}, \ldots, c_{d}) \overset{E_{d+1}(R[X])}{\sim} 
(u_{0}(X), u_{1}(X), c_{2}, \ldots, c_{d}).
 \end{equation} 
 By equation \ref{eq2} and equation \ref{eq3}, we have 
 \begin{equation}
  \label{eq4} 
  (uu_{0}(X), u_{1}(X), c_{2}, \ldots, c_{d}) \overset{E_{d+1}(R[X])}{\sim}  (v_{0}(X), v_{1}(X), \ldots, v_{d}(X)).
\end{equation} 
Therefore by equation \ref{eq1} and equation \ref{eq4}, 
$$(uv_{0}(X), v_{1}(X), \ldots, v_{d}(X))\overset{E_{d+1}(R[X])}{\sim} (v_{0}(X), v_{1}(X), \ldots, v_{d}(X)).$$

$~~~~~~~~~~~~~~~~~~~~~~~~~~~~~~~~~~~~~~~~~~~~~~~~~~~~~~~~~~~~~~~~~~~~~~~~~~~~~~~~~~~~~~~~~~~~~~~~~~~~~~~~~~~~~~\qedwhite$

\begin{notn} Let $R$ be a commutative ring and $T$ be the set of all monic polynomials in $R[X].$ Then we denote 
$R\langle X \rangle = R[X]_{T}$.
\end{notn}

\begin{lem}
 \label{heightlemma} Let $R$ be a commutative noetherian ring of dimension $d$ and $I$ be an ideal of $R[X]$ such 
 that $ht(I) = d.$ Let $T$ be the set of all monic polynomials in $R[X]$ and $J = IR[X]_{T} = IR\langle X\rangle.$ Then 
 $\mbox{ht}(J) \geq d.$
\end{lem}
${\pf}$ Without loss of generality we may assume that $I$ does not contain a monic polynomial. 
\par {\bf{Case I}} : $I$ is a prime ideal. 
\par Since $\mbox{ht}(I) = d,$ there is a chain of prime ideals 
$$\mathfrak{p}_{0} \subsetneq \mathfrak{p}_{1} \subsetneq \cdots \subsetneq \mathfrak{p}_{d} = I.$$ 
Since $I$ does not contain a monic polynomial, we have a chain of prime ideals in 
$R\langle X\rangle,$ 
$$\mathfrak{p}_{0}' \subsetneq \mathfrak{p}_{1}' \subsetneq \cdots \subsetneq \mathfrak{p}_{d}' = J$$ 
where $\mathfrak{p}_{i}' = \mathfrak{p}_{i}R\langle X\rangle.$ Thus $\mbox{ht}(J) \geq d.$ 
\par {\bf{Case II}} : $I$ is not a prime ideal. 
\par Let $\mathfrak{q}$ be a prime ideal containing $IR\langle X\rangle$. Then $\mathfrak{p} = \mathfrak{q}\cap R[X]$ is a prime ideal of $R[X]$ containing $I$, 
therefore its height is $\geq$ $d$ by our assumption on $I.$ By {\bf{Case I}}, the prime ideal $\mathfrak{p}R\langle X\rangle$ has height $\geq$ $d$, it is also contained in $\mathfrak{q}$, which shows that $\mbox{ht}(\mathfrak{q}) \geq d.$ Thus, $\mbox{ht}({IR\langle X\rangle}) \geq d.$

$~~~~~~~~~~~~~~~~~~~~~~~~~~~~~~~~~~~~~~~~~~~~~~~~~~~~~~~~~~~~~~~~~~~~~~~~~~~~~~~~~~~~~~~~~~~~~~~~~~~~~~~~~~~~~~\qedwhite$

Now we prove the Theorem \ref{maintheorem}, stated in the introduction.

{\bf{Proof of Theorem \ref{maintheorem}:}} Let $v = (v_{0}, v_{1}, \ldots, v_{d})\in Um_{d+1}(R[X]).$  In view of $(${\cite[Corollary 1.4]{4}}$),$  $E_{d+1}(R[X]) \trianglelefteq GL_{d+1}(R[X]),$ so neither side 
changes if 
we replace $v$ with $v\alpha$, $g$ with $g\beta\gamma$ where $\alpha, \beta, \gamma \in E_{d+1}(R[X]).$ Apply 
$($\cite[Lemma 3.4]{vdk1}$)$ to $v$ and the first row of $g^{-1}$ to achieve that this first row of $g^{-1}$ equals 
$(w_{0}, v_{1}, \ldots, v_{d}),$ where $w_{0} \in R[X],$ $(v_{1}, \ldots, v_{d})$ is in general position. Denote the first row of $g$ by $(u_{0}, u_{1}, \ldots, u_{d})$ and 
$$ g^{-1} = \begin{bmatrix}
                 w_{0} & v_{1}\cdots v_{d}\\
                         \ast & N\\
                        \end{bmatrix}.$$ 
                    Let $T$ denote the set of all monic polynomials of $R[X].$ Thus $R\langle X\rangle = R[X]_{T}$ and 
                    $\mbox{dim}R\langle X\rangle = d.$ Denote by $v_{i}', u_{i}', w_{i}'$ the images of $v_{i}, u_{i}, w_{i}$ 
                    respectively in $R\langle X\rangle.$ Denote by $N', g'$ the images of $N$ and $g$ respectively in $M_{d}(R\langle X\rangle).$ 
                    \par
                    Let $J = <v_{1}', v_{2}', \cdots, v_{d}'>$, then $\mbox{dim}\frac{R\langle X\rangle}{J} = 0$ by Lemma \ref{heightlemma}. 
                        Let $\mathfrak{m}_{1}, \ldots, \mathfrak{m}_{q}$ be the 
       finitely many maximal ideals of $R\langle X\rangle$ containing $J$ and let $S$ be the complement in $R\langle X\rangle$ of 
       $\bigcup_{i=1}^{q}\mathfrak{m}_{i}.$ 
       Note that $J = \langle u_{1}', u_{2}',\cdots, u_{d}'\rangle $ by using Cramer's rule twice. Let $M'$ be the adjugate matrix of $N'$ 
       such that 
       $N'M' = M'N' = \mbox{det}(N')\cdot I_{d}.$ As $\mbox{det}(g') = a^{-1}$ for some $a\in R^{\times}$, we have 
       $\mbox{det}(N') = au_{0}' \in S$ and $\mbox{det} (M') = (au_{0}')^{d-1} \in S.$ Therefore $M' \in GL_{d}(S^{-1}R\langle X\rangle).$ 
       The ring 
       $S^{-1}R\langle X\rangle$ is semi-local, so 
       $$M' \in \begin{bmatrix}
                 (au_{0}')^{d-1} & 0\\
                         0 & I_{d-1}\\
                        \end{bmatrix}E_{d+1}(S^{-1}R\langle X\rangle).$$ 
         Choose $s\in S$ such that 
       $$M' = \begin{bmatrix}
                 (au_{0}')^{d-1} & 0\\
                         0 & I_{d-1}\\
                        \end{bmatrix}\varepsilon_{1}\cdot \ldots \cdot \varepsilon_{r}\in GL_{d}(R\langle X\rangle[1/s]),$$ 
                         with $\varepsilon_{i}$ elementary matrices. Note that $\overline{s}$ is invertible in 
                        $\frac{R\langle X\rangle}{J}$. Choose $p\in R\langle X \rangle$ with $ps \equiv v_{0}' -w_{0}' ~\mbox{mod}~J.$ As we can add multiples 
                        of $v_{1}', \ldots, v_{d}'$ to $v_{0}'$, we arrange that $ps = v_{0}' - w_{0}'.$ 
                Then 
                \begin{align*}
                 \label{eq5}
                 [v'g'] & = [(w_{0}' + v_{0}' - w_{0}', v_{1}', \ldots, v_{d}')\cdot g']\\
                 & = [1 + (v_{0}' - w_{0}')u_{0}', (v_{0}'-w_{0}')u_{1}', \ldots, (v_{0}'-w_{0}')u_{d}']\\
                 & = [1+psu_{0}', (v_{0}'-w_{0}')u_{0}'^{2}u_{1}', \ldots, (v_{0}'-w_{0}')u_{0}'^{2}u_{d}'].
                \end{align*}
 Now
                \begin{align*}
                [(v_{0}'-w_{0}')&u_{0}'^{2}(u_{1}', \ldots, u_{d}')] \\
                & = [(v_{0}'-w_{0}')u_{0}'(a^{-1}u_{1}', \ldots, a^{-1}u_{d}')N'M']\\
                & = [-(v_{0}'-w_{0}')u_{0}'^{2}(a^{-1}v_{1}', \ldots, a^{-1}v_{d}')M']~~~(\mbox{look~at~the~first~row~of~} g'g'^{-1})\\
                & = [(w_{0}'-v_{0}')u_{0}'^{2}a^{-1}(v_{1}'(au_{0}')^{d-1}, v_{2}', \ldots, v_{d}')]\varepsilon_{1}\cdot \ldots \cdot \varepsilon_{r},
                 \end{align*}
                 where the $\varepsilon_{i}$ are in $E_{d}(R\langle X\rangle)[1/s]),$  not in $E_{d}(R\langle X\rangle).$ On the other hand, we can  
                  multiply coordinates such as $(v_{0}'-w_{0}')u_{0}'^{2}u_{1}'$ freely by $s^{2}$. Therefore, 
                  we mimic the effect of $\varepsilon_{i}$ as follows. For each $\varepsilon_{i}$ choose a diagonal matrix $d_{i}$ with positive 
                  powers of $s^{2}$ on the diagonal, so that $\varepsilon_{i}d_{i} = d_{i}\varepsilon_{i}'$ in $GL_{d}(R\langle X\rangle[1/s])$ for some 
                  $\varepsilon_{i}'\in E_{d}(R\langle X\rangle).$ Choose $d_{i}'$, also a diagonal matrix with positive powers of $s^{2}$ on the 
                  diagonal, so that $d_{i}d_{i}'$ is central in $GL_{d}(R\langle X\rangle[1/s]).$ Then
                  \begin{align*}
                   [(v_{0}'-w_{0}')&u_{0}'^{2}(u_{1}', \ldots, u_{d}')](d_{1}d_{1}'\cdots d_{r}d_{r}') \\
                   & = [(w_{0}'-v_{0}')u_{0}'^{2}a^{-1}(v_{1}'(au_{0}')^{d-1}, v_{2}', \ldots, v_{d}')](d_{1}e_{1}'d_{1}'\cdots 
                   d_{r}e_{r}'d_{r}')
                  \end{align*}
                  (over $R\langle X\rangle[1/s]$ and one can also arrange it to be true over $R\langle X\rangle$ too) so that 
                  $$[v'g'] = [1+psu_{0}', (w_{0}'-v_{0}')u_{0}'^{2}a^{-1}(v_{1}'(au_{0}')^{d-1}, v_{2}', \ldots, v_{d}')].$$ 
                  Thus, 
      \begin{align*}
       [v'g'] & = [1+psu_{0}', (w_{0}'-v_{0}')a^{-1}v_{1}'u_{0}'^{d-1}, (w_{0}'-v_{0}')v_{2}', \ldots, (w_{0}'-v_{0}')v_{d}']\\
           &  = [1+(v_{0}'-w_{0}')u_{0}', ((w_{0}'-v_{0}')u_{0}')^{d-1}a^{-1}(w_{0}'-v_{0}')v_{1}',
          v_{2}', \ldots, v_{d}']\\
          & = [1+(v_{0}'-w_{0}')u_{0}', (w_{0}'-v_{0}')a^{-1}v_{1}',
          v_{2}', \ldots, v_{d}'].
          \end{align*}
          Thus, there exists a monic polynomial $f\in R[X]$ such that
          $$[vg]_{f} = [1+(v_{0}-w_{0})u_{0}, (w_{0}-v_{0})a^{-1}v_{1},
          v_{2}, \ldots, v_{d}]_{f}.$$ 
          Therefore, by $(${\cite[Theorem 1.1]{twoexam}}$)$ 
          $$[vg] = [1+(v_{0}-w_{0})u_{0}, (w_{0}-v_{0})a^{-1}v_{1},
          v_{2}, \ldots, v_{d}].$$ 
          Now by Lemma \ref{removingunit}, 
          $$[vg] = [1+(v_{0}-w_{0})u_{0}, (w_{0}-v_{0})v_{1},
          v_{2}, \ldots, v_{d}].$$ 
         Thus, $[vg] = [v]\ast [-u_{0}, v_{1}, \ldots, v_{d}]~~~~\mbox{as}~~~~(-u_{0})(-w_{0})\equiv 1~\mbox{mod} ~
         \langle v_{1}, v_{2}, \cdots, v_{d}\rangle.$
      \par
      But $[-u_{0}, v_{1}, \ldots, v_{d}] = [w_{0}, v_{1}, \ldots, v_{d}]^{-1} = [e_{1}g^{-1}]^{-1}.$ We have shown that 
      $[vg] = [v]\ast[e_{1}g^{-1}]^{-1}.$ Substituting $(1, 0, \ldots 0) 
      $ for $v$ we see that $[e_{1}g] = [e_{1}g^{-1}]^{-1}$ and the Theorem follows. 
   $~~~~~~~~~~~~~~~~~~~~~~~~~~~~~~~~~~~~~~~~~~~~~~~~~~~~~~~~~~~~~~~~~~~~~~~~~~~~~~~~~~~~~~~~~~~~~~~~~~~~~~~~~~~~~~\qedwhite$
   
   \medskip
\par
   In $(${\cite[Example 4.16]{vdk1}}$),$ W. van der Kallen shows that for a ring $R$ of dimension $d$, the commutator subgroup 
   $[GL_{d+1}(R), GL_{d+1}(R)]$ need not to be contained in $GL_{d}(R)E_{d+1}(R)$. As a consequence of Theorem \ref{maintheorem}, we prove the Corollary \ref{coro}.\\
   {\bf{Proof of Corollary \ref{coro}:}} In view of Theorem \ref{maintheorem}, the first row map 
 $$GL_{d+1}(R[X]) \longrightarrow \frac{Um_{d+1}(R[X])}{E_{d+1}(R[X])}$$  
 is a group homomorphism. Thus, we have $$[GL_{d+1}(R[X]), GL_{d+1}(R[X])] \subseteq GL_{d}(R[X])E_{d+1}(R[X]).$$ 
 $~~~~~~~~~~~~~~~~~~~~~~~~~~~~~~~~~~~~~~~~~~~~~~~~~~~~~~~~~~~~~~~~~~~~~~~~~~~~~~~~~~~~~~~~~~~~~~~~~~~~~~~~~~~~~~\qedwhite$

  \medskip
\noindent
{\bf Acknowledgement:} I thank the referee for indicating 
an easy proof of Lemma \ref{heightlemma} and for going through
the manuscript with great care. A detailed list of suggestions by the referee improved
the exposition considerably.

\Addresses

\end{document}